\documentclass{amsart}
\usepackage{amsmath,amssymb}
\usepackage[all]{xy}

\theoremstyle{plain}
\newtheorem{theorem}{Theorem}[section]
\newtheorem{proposition}[theorem]{Proposition}
\newtheorem{lemma}[theorem]{Lemma}
\newtheorem{corollary}[theorem]{Corollary}

\theoremstyle{definition}
\newtheorem{definition}[theorem]{Definition}
\newtheorem{example}[theorem]{Example}

\theoremstyle{remark}

\newcommand{\secref}[1]{Section~\ref{#1}}
\newcommand{\thmref}[1]{Theorem~\ref{#1}}
\newcommand{\propref}[1]{Proposition~\ref{#1}}
\newcommand{\lemref}[1]{Lemma~\ref{#1}}
\newcommand{\corref}[1]{Corollary~\ref{#1}}

\newcommand{\defref}[1]{Definition~\ref{#1}}

\def\R{{\mathbb \R}}

\def\Q{{\mathbb Q}}
\def\C{{\mathbb C}}

\def\map{\mathrm{map}}

\def\Hom{\mathrm{Hom}}

\def\cat0{\mathrm{cat}_0}

\begin{document}

\title[Cyclic Maps in Rational Homotopy Theory]{Cyclic Maps in Rational
Homotopy Theory}

\author{Gregory  Lupton}

\address{Department of Mathematics,
          Cleveland State University,
          Cleveland OH 44115}

\email{G.Lupton@csuohio.edu}

\author{Samuel Bruce Smith}

\address{Department of Mathematics,
  Saint Joseph's University,
  Philadelphia, PA 19131}

\email{smith@sju.edu}

\date{\today}



\begin{abstract}
The notion of a cyclic map $g \colon  A \to X$ is a natural
generalization of a Gottlieb element in $\pi_n(X)$. We investigate
cyclic maps from a rational homotopy theory point of view. We show
a number of results for rationalized cyclic maps which generalize
well-known results on the rationalized Gottlieb groups.
\end{abstract}

\maketitle

\section{Introduction}\label{sec:intro}

Let $X$ be a space with basepoint $x_0$.  Denote by $\map(X, X;1)$
the component of the identity in the unbased function space of
maps from $X$ to itself, and by $\omega\colon  \map(X, X;1)
\rightarrow X$ the evaluation map defined by $\omega(f)(x) =
f(x_0)$. In \cite{Go1}, Gottlieb introduced and studied the
evaluation subgroups
$$G_n(X) = \omega_\#\big(\pi_n\big(\map(X, X;1)\big)\big)\subseteq \pi_n(X).$$
Note that $G_n(X)$ can alternatively be described as homotopy
classes of maps $f \colon  S^n \to X$ such that $(f \mid 1) \colon
S^n \vee X \to X$ admits an extension $F \colon S^n\times X \to X$
up to homotopy.  These so-called \emph{Gottlieb groups} admit a
number of generalizations.  One such was introduced and studied by
Varadarajan.  In \cite{Var}, he defines the homotopy set $G(A, X)$
of \emph{cyclic maps} from $A$ to $X$, that is, homotopy classes
of maps $f \colon A \to X$ such that $(f \mid 1) \colon A \vee X
\to X$ admits an extension $F \colon A \times X \to X$. The
Gottlieb group occurs when $A=S^n$.  Note that the set of cyclic
maps $G(A, X)$ can alternatively be described as the set of
homotopy classes of maps from $A$ to $X$ that admit a lift through
the evaluation map $\omega\colon  \map(X, X;1) \rightarrow X$.

The rationalized Gottlieb groups $G_n(X)\otimes\Q$ have been
studied extensively (see \cite{F-H, Tan83, Oprea}). Our purpose in
this paper is to study rationalized cyclic maps with an eye toward
obtaining natural generalizations of known results. In \cite{Go1}
Gottlieb studied the question of when the Gottlieb group of $X$ is
contained in the kernel of the mod $p\,$  or rational Hurewicz
homomorphisms. Let $h_\infty \colon \pi_*(X) \to H_*(X; \Q)$
denote the rational Hurewicz homomorphism.  Then Gottlieb's main
results for the rational case can be summarized as follows:

\begin{theorem}[{\cite[Th.4.1,Th.5.1]{Go1}}]\label{thm:Gottlieb}%
Let $X$ be a space with finitely generated homology.  Then
$G_{2n}(X) \subseteq \ker ( h_\infty )$ and if the Euler
characteristic $\chi(X) \neq 0$ then $G_{2n-1}(X) \subseteq
\ker(h_\infty)$.
\end{theorem}

Gottlieb also showed that $G_{2}(X)$ is a torsion group---and
hence that $G_{2}(X)\otimes\Q = 0$---in the case that $X$ has a
finite number of non-zero rational homology groups
\cite[Th.7.1]{Go1}.

F{\'e}lix and Halperin significantly sharpened these results. Let
$X$ be a simply connected space.  Then the \emph{rational
category} of $X$---denoted by $\cat0(X)$---may be defined as the
Lusternik-Schnirelmann category of the rationalization of $X$,
that is, $\cat0(X) = \mathrm{cat}(X_\Q)$. Further, the
rationalized Gottlieb group of $X$ is contained in the Gottlieb
group of the rationalization of $X$, or $G_n(X)\otimes\Q \subseteq
G_n(X_\Q)$ (these agree if $X$ is finite \cite{Lan}).

\begin{theorem}[{\cite[Th.III]{F-H}}]\label{thm:Felix-Halperin}%
Let $X$ be a simply connected space of finite rational category.
Then $G_{2n}(X_\Q) = 0$ and $\sum_n
\mathrm{dim}\big(G_{2n+1}(X_\Q)\big) \leq \mathrm{cat}_0(X)$.
\end{theorem}

In this paper, our main results are presented in
\secref{sec:factor}. There, we give a generalization to cyclic
maps of the F{\'e}lix-Halperin result on even dimensional Gottlieb
groups. \corref{cor:rational cyclic} shows that cyclic maps of
rational spaces are trivial under certain conditions.  That result
includes the part of \thmref{thm:Felix-Halperin} about trivial
even dimensional Gottlieb groups of a rational space as a special
case.

We also generalize Gottlieb's result about odd dimensional
Gottlieb groups, as follows: To say that $\alpha \in \pi_n(X)$ is
in the kernel of the rational Hurewicz homomorphism is simply to
say that $H_*(\alpha) = 0 \colon H_*(S^n; \Q) \to H_*(X; \Q)$.
More generally, write $SH_*(X; \Q)$ for the  {\em rational
spherical homology} of $X$: that is, $SH_*(X; \Q)$ is   the image
of the rationalized Hurewicz homomorphism $h_\Q : \pi_*(X) \otimes
\Q \to H_*(X; \Q).$
  Note that a map $f : A \to X$ induces a map
$SH_*(f) : SH_*(A; \Q) \to SH_*(X; \Q).$
  We
give a vanishing result for the map induced on rational spherical homology
by a cyclic map in \corref{cor:chi non-zero}, which extends
Gottlieb's  odd dimensional result.

Also in \secref{sec:factor}, we show results about the
homomorphism induced on rational homotopy groups by a cyclic map.
One such result is \corref{cor:rank f sharp}, which gives a bound
on the size of the image that a cyclic map may induce on odd
dimensional rational homotopy groups. It complements the
F{\'e}lix-Halperin result on odd dimensional Gottlieb groups. Many
of the results in \secref{sec:factor} are actually proved for a
more general class of maps than cyclic maps, namely, maps that
factor through an $H$-space.

In \secref{sec:examples} we give a number of examples of cyclic
maps between rational spaces to
illustrate our results in \secref{sec:factor}. We also give some
sample  computations of the rationalization of  the set of cyclic
maps.

All spaces in the sequel are assumed to be simply connected CW
complexes with rational homology of finite type---with the
exception of mapping spaces such as $\map(X,X;1_X)$. We denote the
set of homotopy classes of maps from $X$ to $Y$ by $[X,Y]$. We
often do not distinguish between a map and the homotopy class it
represents. If $f \colon A \to B$ is a map, then $f^*$ denotes
pre-composition by $f$ and $f_*$ denotes post-composition by $f$.
Thus we obtain maps of homotopy sets $f^*\colon [B,X] \to [A,X]$
and $f_*\colon [X,A] \to [X,B]$. We use $H_*(f)$ and $H^*(f)$ to
denote the map induced on homology, respectively cohomology, by
the map of spaces $f$, and $f_\#$ to denote the map induced on
homotopy groups.   We denote the rationalization of a space $X$ by
$X_\Q$ and of a map $f$ by $f_\Q$ (cf. \cite{H-M-R}).

We assume some familiarity with rational homotopy theory as
introduced by Sullivan. Our main reference for this material is
\cite{F-H-T}. The basic facts that we use are as follows: Each
space $X$ has a unique Sullivan minimal model $(\mathcal{M}_{X},
d_X)$ in the category of simply connected DG (differential graded)
algebras over $\Q$.  This DG algebra $(\mathcal{M}_{X}, d_X)$ is
of the form $\mathcal{M}_{X} = \Lambda V$, a free graded
commutative algebra generated by a positively graded vector space
$V$ of finite type.  The differential $d_X$ is decomposable, in
that $d_X(V) \subseteq \Lambda^{\geq 2} V$. A map $f\colon X \to
Y$ has a Sullivan minimal model which is a DG algebra map
$\mathcal{M}_f \colon \mathcal{M}_Y \to \mathcal{M}_X$. The
Sullivan minimal model is a complete rational homotopy invariant
for a space or a map. Since the minimal model is determined by the
rational homotopy type, the minimal models of $X_\Q$ and $X$, and
more generally those of $f_\Q$ and $f$, agree. There is a notion
of homotopy for maps of Sullivan minimal models, which we refer to
as \emph{DG homotopy} (of maps of Sullivan minimal models). If $f,
g \colon X \to Y$ are maps of \emph{rational spaces}, then $f$ and
$g$ are homotopic if and only if their Sullivan minimal models
$\mathcal{M}_f$ and $\mathcal{M}_g$ are DG homotopic. Rational
cohomology is readily retrieved from Sullivan minimal models:  We
have a natural isomorphism $H(\mathcal{M}_{X}, d_X) \cong
H^*(X;\Q)$ and this isomorphism identifies $H(\mathcal{M}_f)
\colon H(\mathcal{M}_Y) \to H(\mathcal{M}_X)$ with $H^*(f) \colon
H^*(Y; \Q) \to H^*(X; \Q)$. Rational homotopy is retrieved as
follows:  Let $Q(\mathcal{M}_X) \cong V$ be the (quotient) module
of indecomposables of $\mathcal{M}_X$. There is a natural
isomorphism $Q(\mathcal{M}_X) \cong \Hom(\pi_*(X)\otimes\Q,\Q)$,
that identifies $Q(\mathcal{M}_f) \colon Q(\mathcal{M}_Y) \to
Q(\mathcal{M}_X)$ with $(f_\#\otimes\Q)^* \colon
\Hom(\pi_*(Y)\otimes\Q,\Q) \to \Hom(\pi_*(X)\otimes\Q,\Q)$.

For facts that we have used about minimal models beyond these
basics, we have given specific references from \cite{F-H-T}.   We
use the standard notation and terminology for minimal models as
used in \cite{F-H-T}. A map of DG algebras is called a
\emph{quasi-isomorphism} if it induces an isomorphism on
cohomology. We mention that there is a discussion of Gottlieb
groups, including a proof of \thmref{thm:Felix-Halperin} cited
above, given in \cite[sec.28(d)]{F-H-T}.

\section{Generalities on Cyclic Maps}

In this section,  we record some basic facts about the set
$G(A,X)$. All of these results are either known as stated or occur
as easy extensions of known results on cyclic maps or evaluation
subgroups. The first is a useful result from \cite{Var}.

\begin{theorem}[{\cite[Th.1.3]{Var}}]\label{thm:cyclic composition}%
Let $ \theta \colon  B \to A$ be any  map. Then $\theta^* \colon
[A,X] \to [B,X]$ restricts to a map of sets of cyclic maps
$\theta^* \colon G(A,X) \to G(B,X)$.
\end{theorem}

\begin{proof} Suppose that $F \colon  A \times X \to X$ extends
$(g \mid 1)$.  Then $G\colon B \times X \to X$ defined by $G(b, x)
= F(\theta(b), x)$ extends $(g \circ \theta \mid 1)$.
\end{proof}

We obtain the following two immediate consequences of this.

\begin{corollary}\label{cor:image of cyclic is Gottlieb}
A cyclic map $g  \in G(A,X)$ satisfies
$$g_\sharp(\pi_n(A)) \subseteq G_n(X).$$
\end{corollary}

\begin{proof} Take $\theta \colon  S^n \to A$. Then  $g_\#(\theta) = g \circ
\theta = \theta^*(g)\colon  S^n \to X$ is in $G_n(X)$.
\end{proof}

\begin{corollary} $G(X, X)$ contains  a self-equivalence
of $X$ if and only if $X$ is an $H$-space.
\end{corollary}

\begin{proof} Let $g \colon  X \to X$ be a cyclic self-equivalence and $h \colon
X \to X$ its homotopy inverse.  Then, by Theorem 2.1, $1_X \simeq
g \circ h$ is cyclic as well.  The map $F\colon  X \times X \to X$
that extends $1_X \vee 1_X$ is the needed  multiplication.
\end{proof}

Next we give some results, and establish our notation, concerning
localization (see \cite{H-M-R} for full details). For any set of
primes $P$, the $P$-localization map $e_X \colon X \to X_P$
induces a map of homotopy sets $(e_X)_* \colon [A,X] \to [A,X_P]$.
Likewise, we have an induced map $(e_A)^* \colon [A_P,X_P] \to
[A,X_P]$, which is a bijection of sets.

\begin{lemma}\label{lem:localization}
Let $A$ and $X$ be simply connected and let $P$ be any set of
primes.
\begin{enumerate}
\item $(e_X)_*$ restricts to a map $(e_X)_* \colon G(A,X) \to
G(A,X_P)$.
\item The bijection $(e_A)^*$ restricts to a bijection $(e_A)^* \colon G(A_P,X_P) \to
G(A,X_P)$.
\end{enumerate}
\end{lemma}

\begin{proof}
Both parts follow easily from standard properties of localization.
Suppose that $f \colon A \to X$ is a cyclic map with affiliated
map $F \colon A \times X \to X$.  Then an affiliated map for
$(e_X)_*(f)$ is given by $F_P \circ (e_A \times 1) \colon A \times
X_P \to X_P$.  This shows (1).  For (2), suppose $\alpha \colon
A_P \to X_P$ is a cyclic map. Then $(e_A)^*(\alpha) \in G(A,X_P)$
by \thmref{thm:cyclic composition}. Since the restriction of an
injection is injective, it only remains to show that $(e_A)^*$
restricts to a surjection.  For this, suppose $\beta \in G(A,X_P)$
has affiliated map $B \colon A \times X_P \to X_P$.  Then $\beta =
(e_A)^*(\beta_P)$, and $B_P \colon A_P \times X_P \to X_P$ is an
affiliated map for $\beta_P$. The latter assertion uses uniqueness
of localized maps, and the fact that both $e_A \vee 1 \colon A
\vee X_P \to A_P \vee X_P$ and $e_A \times 1 \colon A \times X_P
\to A_P \times X_P$ are $P$-localization maps.
\end{proof}

Because we use minimal model techniques, we obtain results about
homotopy classes of cyclic maps $G(A_\Q,X_\Q)$, or more generally
about homotopy classes of maps $[A_\Q,X_\Q]$.

\begin{definition}\label{def:rationally cyclic}
Let $f \colon A \to X$ be a map of simply connected spaces.  We
say $f$ is a \emph{rationally cyclic map} if its rationalization
$f_\Q \colon A_\Q \to X_\Q$ is a cyclic map. More generally, we
say $f$ is a $P$-local cyclic map if $f_P \colon A_P \to X_P$ is
cyclic. We denote the set of homotopy classes of $P$-local cyclic
maps from $A$ to $X$ by $G_P(A,X)$.  From
\lemref{lem:localization}, we have $G(A,X) \subseteq G_P(A,X)$ for
any set of primes $P$.
\end{definition}

In general the inclusion $(e_X)_*\big(G(A,X)\big) \subseteq
G(A,X_P)$ may be strict: See \cite[pp.378--380]{F-H-T} for an
(infinite) example with $(e_X)_*\big(G_n(X)\big) \subsetneq
G_n(X_\Q)$, but note that $(e_X)_*\big(G_n(X)\big) = G_n(X_P)$ if
$X$ is finite, by \cite{Lan}.  We do have the following result,
however.

\begin{lemma}\label{lem:finiteness}
Let $A$ and $X$ be simply connected with $A$ a finite CW complex.
Then $P$-localization induces a finite-to-one map
$$\big((e_A)^*\big)^{-1}\circ (e_X)_* \colon G_P(A,X) \to G(A_P,X_P).$$
In particular, if $A$ is finite and $G(A_\Q,X_\Q)$ is trivial,
that is, consists of a single element, then the set of rationally
cyclic maps $G_\Q(A,X)$ is finite.
\end{lemma}

\begin{proof}
The finite-to-one assertion follows from \cite[Cor.5.4]{H-M-R}.
\end{proof}

\section{Maps that Factor Through an $H$-space}\label{sec:factor}

In this section we show that being cyclic entails strong
restrictions on a map. Most of our results flow from the following
observation. Suppose that $f \colon A \to X$ is a cyclic map.  As
we remarked in the introduction, $f$ admits a lift through the
evaluation map as follows
$$\xymatrix{ & \map(X,X;1) \ar[d]^{\omega} \\
A \ar[ur]^{\tilde f} \ar[r]_{f}& X.}$$
Now the identity component $\map(X,X;1)$ is an $H$-space, and thus
a cyclic map factors through an $H$-space.  Indeed, since we are
assuming $A$ is simply connected, we can lift $\tilde f$ through
the universal cover of $\map(X,X;1)$, which is again an $H$-space,
and is also simply connected.  Thus \emph{a cyclic map factors
through a simply connected $H$-space}.  Such a factorization
entails numerous consequences for a cyclic map. In fact, some of
our results can be proved simply by assuming that the map $f
\colon A \to X$ factors through some $H$-space---not necessarily
the universal cover of $\map(X,X;1)$. This hypothesis on a map is
definitely weaker than assuming the map to be cyclic. For
instance, consider any map $f \colon S^3 \to S^3 \vee S^3$ that is
not trivial on homotopy. Since $G_3(S^3 \vee S^3) = 0$, such a map
is not a cyclic map. Since $S^3$ itself is an $H$-space, any such
map certainly factors through an $H$-space.

Our primary interest is in drawing conclusions about cyclic maps.
However, since our methods are those of rational homotopy, we need
only require that the map be a rationally cyclic map as defined in
\defref{def:rationally cyclic}. Likewise, although a cyclic map
actually factors through an $H$-space, we only require such a
factorization after rationalization. In what follows, we will
state our results so as to place as weak a hypothesis as possible
on the maps and spaces.  We remind the reader that our spaces are
assumed to be simply connected with rational homology of finite
type.

Our first consequence of these observations is as follows.

\begin{proposition}\label{prop:image in cycles}
Let $f \colon A \to X$ be a map whose rationalization factors
through an $H$-space. Up to DG homotopy, the Sullivan minimal
model $\mathcal{M}_f \colon \mathcal{M}_X \to \mathcal{M}_A$ of
$f$ has image contained in the cycles of $\mathcal{M}_A$.
\end{proposition}

\begin{proof}
Suppose we have a factorization
\begin{equation}\label{eqn:factor}
\xymatrix{ & Y \ar[d]^{g} \\
A_\Q \ar[ru]^{\tilde f} \ar[r]_{f_\Q} & X_\Q}
\end{equation}
with $Y$ an $H$-space. Passing to minimal models, we have
$\mathcal{M}_f \sim \mathcal{M}_{\tilde f}\circ \mathcal{M}_g
\colon \mathcal{M}_{X} \to \mathcal{M}_{A}$.  Now the minimal
model of an $H$-space has trivial differential, that is, $d_Y = 0$
\cite[p.143]{F-H-T}.  For every element $\chi \in
\mathcal{M}_{X}$, therefore, we have that $\mathcal{M}_{g}(\chi)$,
and hence $\mathcal{M}_{\tilde f}\circ \mathcal{M}_g(\chi)$, is a
cycle.
\end{proof}

\propref{prop:image in cycles} gives a serendipitious
justification of the terminology ``cyclic map,"  at least in
rational homotopy (cf.~the remarks on nomenclature in \cite{Var}).
With additional hypotheses on either $A$ or $X$, we can draw much
stronger conclusions.

\begin{theorem}\label{thm:odd rationally trivial}
Let $f \colon A \to X$ be a map whose rationalization factors
through an $H$-space. If $X$ is of finite rational category and
$H^{\mathrm{odd}}(A;\Q) = 0$, then $f$ is rationally trivial.
\end{theorem}

\begin{proof}
We prove the result by showing that the minimal model
$\mathcal{M}_f\colon \mathcal{M}_X \to \mathcal{M}_A$ of $f$ is DG
homotopic to the trivial map.  See \cite[Sec.12(b)]{F-H-T} for
details about DG homotopy.

Suppose we have a factorization as in (\ref{eqn:factor}), with $Y$
an $H$-space.  First we show that, up to DG homotopy, we can
assume that $\mathcal{M}_{\tilde f}$ is zero on all odd-degree
generators of $\mathcal{M}_{Y}$.  For suppose that
$\mathcal{M}_{Y} = \Lambda(\{a_i\}_{i\in I}, \{b_j\}_{j\in J})$,
with each $a_i$ an even-degree generator and each $b_j$ an
odd-degree generator.  Recall that the differential $d_Y$ is zero,
and so every element in $\mathcal{M}_{Y}$ is a cycle. Since
$\tilde f(b_i)$ is an odd-degree cycle in $\mathcal{M}_{A}$, and
since we are assuming that $H^{\mathrm{odd}}(A;\Q) = 0$, we have
$\tilde f(b_j) = d_A(\eta_j)$ for some $\eta_j \in
\mathcal{M}_{A}$.  Now define a map $\Phi \colon \mathcal{M}_{Y}
\to \mathcal{M}_{A}\otimes\Lambda(t,dt)$ on generators by setting
$\Phi(b_j) = \tilde f(b_j)\otimes(1-t) + \eta_j\otimes dt$,
$\Phi(a_i) = \tilde f(a_i)\otimes1$.  As the differential in
$\mathcal{M}_{Y}$ is trivial, this map extends to a DG homotopy
that starts at $\tilde f$ and ends at a map that is zero on each
odd-degree generator $b_j$.  So now assume that
$\mathcal{M}_{\tilde f}$ is zero on odd-degree generators.

Let $\mathcal{I}$ be the ideal of $\mathcal{M}_{Y}$ generated by
the odd-degree generators.  We now show that $\mathcal{M}_{g}$ has
image contained in $\mathcal{I}$. Clearly, for any odd-degree
generator $x$ of $\mathcal{M}_{X}$, we must have
$\mathcal{M}_{g}(x) \in \mathcal{I}$.  To show that
$\mathcal{M}_{g}(x) \in \mathcal{I}$ for $x$ an even-degree
generator of $\mathcal{M}_{X}$, we argue inductively over degree.
Assume that $\mathcal{M}_{g}(y) \in \mathcal{I}$, for all
generators $y$ of degree $\leq 2n-1$, and that $x$ is of degree
$2n$.  Since $X$ has finite rational category, there exists an
element $\eta \in \mathcal{M}_X$ such that $d_X(\eta) = x^k + R$,
for some $k \geq 2$ and some $R$ in the ideal of $\mathcal{M}_{X}$
generated by generators of degree $\leq 2n-1$. We justify this
assertion in \lemref{lem:technical} below. Now write
$\mathcal{M}_{g}(x) = P + Q$, with $Q \in \mathcal{I}$ and $P$ an
element of the subalgebra of $\mathcal{M}_{Y}$ generated by the
even-degree generators. Since the differential in
$\mathcal{M}_{Y}$ is trivial, we have
$\mathcal{M}_{g}\big(d_X(\eta)\big) =
d_Y\big(\mathcal{M}_{g}(\eta)\big) = 0$.  Hence, $0 =
\mathcal{M}_{g}(x^k + R) = \big(\mathcal{M}_{g}(x)\big)^k +
\mathcal{M}_{g}(R) = (P+ Q)^k + \mathcal{M}_{g}(R)$.  Since $Q \in
\mathcal{I}$, and by our inductive assumption $\mathcal{M}_{g}(R)
\in \mathcal{I}$, it follows that $P^k \in \mathcal{I}$ and hence
that $P = 0$.  This completes the inductive step. We can start the
induction with $n = 1$, where the induction hypothesis is
satisfied trivially.  We have shown that $\mathcal{M}_{g}$ has
image contained in $\mathcal{I}$. Since $\tilde f$ is zero on
odd-degree generators, we have $\mathcal{M}_{\tilde
f}\circ\mathcal{M}_{g} = 0$.
\end{proof}

The technical fact about minimal models used in the preceding
proof is well-known amongst experts in the field. We give a
statement and proof of this useful fact for the sake of
completeness.

\begin{lemma}\label{lem:technical}
Let $X$ be a space of finite rational category. If $x \in
\mathcal{M}_X$ is a generator of even degree $2n$, then for some
$k \geq 2$ there is an element $\eta \in \mathcal{M}_X$ such that
$d_X(\eta) = x^k + R$, with $R$ in the ideal of $\mathcal{M}_{X}$
generated by generators of degree $\leq 2n-1$.
\end{lemma}

\begin{proof}
First notice that if the even-degree generator is a cycle, $d_X(x)
= 0$, then the result is easily proved.  For the assumption of
finite rational category implies, in particular, that $X$ has
finite rational cup length.  Hence in this case, there is some
$\eta$ such that $d_X(\eta) = x^k$ with $k \geq 2$. For the
general case, in which $x$ is not a cycle, we use the mapping
theorem of F{\'e}lix-Halperin \cite[Th.29.5]{F-H-T}. Suppose $x$
is of degree $2n$, and let $\mathcal{M}_X\langle 2n \rangle$
denote the quotient DG algebra obtained by factoring out the DG
ideal of $\mathcal{M}_X$ generated by all generators of degree
$\leq 2n-1$. Let $\overline{d_X}$ denote the differential induced
on the quotient by $d_X$.  Then $(\mathcal{M}_X\langle 2n \rangle,
\overline{d_X})$ is a minimal DG algebra and the projection
$\mathcal{M}_X \to \mathcal{M}_X\langle 2n \rangle$ induces a
surjection of the modules of indecomposables.  It follows from the
mapping theorem that $(\mathcal{M}_X\langle 2n \rangle,
\overline{d_X})$ also has finite rational category.  By the
argument at the start of this proof, there is some element $\eta
\in \mathcal{M}_X\langle 2n \rangle$ such that
$\overline{d_X}(\eta) = x^k$ in $\mathcal{M}_X\langle 2n \rangle$.
Therefore, $d_X(\eta) = x^k + R$ for some $R$ in the ideal of
$\mathcal{M}_{X}$ generated by generators of degree $\leq 2n-1$,
as asserted.
\end{proof}

\begin{corollary}[to \thmref{thm:odd rationally trivial}]\label{cor:rational cyclic}
If $\cat0(X)<\infty$ and $H^{\mathrm{odd}}(A;\Q) = 0$, then
$G(A_\Q,X_\Q)$ is trivial.  If, further, $A$ is a finite CW
complex, then $G_\Q(A,X)$, and hence $G(A,X)$, is a finite set.
\end{corollary}

\begin{proof}
As we observed above, any cyclic map---including those in
$G(A_\Q,X_\Q)$---factors through an $H$-space.  Now apply
\thmref{thm:odd rationally trivial} and \lemref{lem:finiteness}.
\end{proof}

\corref{cor:rational cyclic} is a most satisfactory generalization
to cyclic maps of the F{\'e}lix-Halperin result about even
dimensional rational Gottlieb groups.  Unfortunately, we have not
found such a satisfactory generalization of their odd dimensional
result. Our next three results take a step in this direction,
however, giving restrictions on the homomorphism induced on
homotopy groups by a cyclic map.

\begin{proposition}[cf.~{\cite[Cor.~on p.379]{F-H-T}}]\label{prop:rational Gottlieb}
Let $f \colon A \to X$ be a rationally cyclic map, with $X$ a
space of finite rational category.  Then $f_\#\otimes\Q$ is zero
in all even-degrees and $\mathrm{rank}(f_\#\otimes\Q) \leq
\mathrm{cat}_0(X)$.
\end{proposition}

\begin{proof}
By assumption, $f_\Q \colon A_\Q \to X_\Q$ is a cyclic map. Hence
$(f_\Q)_\# \colon \pi_n(A_\Q) \to \pi_n(X_\Q)$ has image in
$G_n(A_\Q)$ by \lemref{cor:image of cyclic is Gottlieb}. Since
$(f_\Q)_\#$ is identified with $f_\# \otimes\Q$, both assertions
now follow by \cite[Th.III]{F-H}.
\end{proof}

We also have the following result that extends the first assertion
of \propref{prop:rational Gottlieb}, and gives a further
restriction on the rank of the homomorphism induced on rational
homotopy groups by a cyclic map.

Recall that we defined the rational spherical homology of $X$ to
be $SH_{n}(X;\Q)= h_\Q (\pi_{n}(X)\otimes\Q)$ for each $n$, where
$h_\Q \colon \pi_*(X)\otimes\Q \to H_*(X;\Q)$ is the rationalized
Hurewicz homomorphism.    The vector space dual of $h_\Q$ gives a
homomorphism $(h_\Q)^* \colon H^*(X;\Q) \to
\Hom(\pi_*(X)\otimes\Q,\Q)$.  The subspace of $H^*(X;\Q)$ dual to
$SH_{*}(X;\Q)$ is referred to as the \emph{rational spherical
cohomology} of $X$, and is denoted by $SH^*(X;\Q)$.  Notice that
$(h_\Q)^*$ restricts to an injection from $SH^{n}(X;\Q)$ into
$\Hom(\pi_{n}(X)\otimes\Q,\Q)$.  Now suppose that $(\mathcal{M}_X,
d_X)$ is the minimal model of $X$.  Since the differential $d_X$
is decomposable, it induces the trivial differential on the module
of indecomposables $Q(\mathcal{M}_X)$.  Therefore, by passing to
cohomology from the quotient projection $\mathcal{M}_X \to
Q(\mathcal{M}_X)$, we obtain a map $\zeta_X \colon
H(\mathcal{M}_X) \to Q(\mathcal{M}_X)$.  Under the natural
identifications of $H(\mathcal{M}_X)$ with $H^*(X;\Q)$ and
$Q(\mathcal{M}_X)$ with $\Hom(\pi_*(X)\otimes\Q,\Q)$, the map
$\zeta_X$ is naturally identified with the dual of the
rationalized Hurewicz homomorphism $(h_\Q)^*$ \cite[p.173]{F-H-T}.

\begin{theorem}\label{thm:cyclic on rational homotopy}%
Let $f \colon A \to X$ be a map whose rationalization factors
through an $H$-space.
\begin{enumerate}
\item $f_\#\otimes\Q$ is zero on $\mathrm{kernel}\big(h_\Q \colon
\pi_{n}(A)\otimes\Q \to H_{n}(A;\Q)\big)$, for each $n$.
\item If $X$ is a space of finite rational category, then
$f_\#\otimes\Q$ is zero in all even-degrees.
\end{enumerate}
\end{theorem}

\begin{proof}
Recall that $f_\#\otimes\Q$ is identified with the map induced by
$\mathcal{M}_f \colon \mathcal{M}_X \to \mathcal{M}_A$ on the
(quotient) modules of indecomposables.  We observed in
\propref{prop:image in cycles} that the minimal model
$\mathcal{M}_f \colon \mathcal{M}_X \to \mathcal{M}_A$ has image
in the cycles of $\mathcal{M}_A$. This means, in particular, that
any indecomposable terms occurring in the image of $\mathcal{M}_f$
must be indecomposable cycles in $\mathcal{M}_A$. Now the vector
space of indecomposable  cycles in $\mathcal{M}_A$ is isomorphic
to the rational spherical cohomology of $A$. Assertion (1)
follows.

Now suppose that $\cat0(X)<\infty$ and that, as in the proof of
\thmref{thm:odd rationally trivial}, we have a factorization
$\mathcal{M}_f \sim \mathcal{M}_{\tilde f}\circ\mathcal{M}_g
\colon \mathcal{M}_X \to \mathcal{M}_A$, for $\mathcal{M}_g \colon
\mathcal{M}_X \to \mathcal{M}_Y$ with $Y$ an $H$-space. In that
proof, we showed that $\mathcal{M}_{g}$ has image contained in
$\mathcal{I}$, the ideal of $\mathcal{M}_{Y}$ generated by the
odd-degree generators (we did not use any hypothesis on $A$ for
that part of the argument). In particular, if $x$ is an
even-degree generator of $\mathcal{M}_{X}$, then
$\mathcal{M}_{g}(x)$ is decomposable in $\mathcal{M}_{Y}$.  It
follows that $g$, and hence $f$, induces zero on all even-degree
rational homotopy groups.  This shows assertion (2).
\end{proof}

Thus we can add the following to \propref{prop:rational Gottlieb}:

\begin{corollary}\label{cor:rank f sharp}
Let $f \colon A \to X$ be a rationally cyclic map with $X$ a space
of finite rational category. Then $\mathrm{rank}(f_\#\otimes\Q)
\leq
\sum_{n \text{ odd}}
\mathrm{dim}\big(SH_{n}(A;\Q)\big)$.
\end{corollary}

Finally, we turn to our  generalization of Gottlieb's odd-degree result.
We first observe a restriction on    the homomorphism induced on rational
cohomology by maps which factor through an $H$-space.

\begin{theorem}\label{thm:induced on even cohomology}
Suppose $X$ has finite rational category.  If $f \colon A \to X$
is a map whose rationalization factors through an $H$-space, then
$H^*(f)\big(H^{\text{even}}(X;\Q)\big) \subseteq H^{+}(A;\Q)\cdot
H^{+}(A;\Q)$.
\end{theorem}

\begin{proof}
We show the minimal model counterpart of the assertion. It follows
directly from \propref{prop:image in cycles} that the image of any
cohomology class of $H^*(\mathcal{M}_X)$ that is represented by a
decomposable cycle in $\mathcal{M}_X$ is contained in
$H^{+}(\mathcal{M}_A)\cdot H^{+}(\mathcal{M}_A)$.  To handle the
case in which a cohomology class is represented by an
indecomposable cycle, we use a fact established in the proof of
\thmref{thm:odd rationally trivial}. As in that proof, suppose we
have a factorization $\mathcal{M}_f \sim \mathcal{M}_{\tilde
f}\circ\mathcal{M}_g \colon \mathcal{M}_X \to \mathcal{M}_A$,
where $\mathcal{M}_g \colon \mathcal{M}_X \to \mathcal{M}_Y$ for
$Y$ an $H$-space.  Let $\mathcal{I}$ denote the ideal of
$\mathcal{M}_{Y}$ generated by the odd-degree generators. Then we
showed in \thmref{thm:odd rationally trivial} that
$\mathcal{M}_{g}$ has image contained in $\mathcal{I}$.  From this
it follows that the image under $\mathcal{M}_{g}$ of any even
degree generator is decomposable in the cycles of $\mathcal{M}_A$.
Consequently, we have
$H(\mathcal{M}_{f})\big(H^{\text{even}}(\mathcal{M}_X)\big)
\subseteq H^{+}(\mathcal{M}_A)\cdot H^{+}(\mathcal{M}_A)$.
\end{proof}

We do not obtain a restriction comparable to that of
\thmref{thm:induced on even cohomology} for odd-degree cohomology.
This is because any cohomology class of $H^{2n+1}(A;\Q)$ is
in the image of a homomorphism induced by a map $f \colon A \to
K(\Q, 2n+1)$, and such a map is a cyclic map.

In the following result, we do not assume that $X$ has finite
rational category. Recall that $SH_*(f) : SH_*(A, \Q) \to SH_*(X,
\Q)$ denotes  the map induced on rational spherical homology by $f : A \to
X.$

\begin{theorem}\label{thm:spherical on cohomology}
Let $f\colon A \to X$ be a rationally cyclic map. If $SH_*(f) \neq  0$
then $X$ decomposes up to rational homotopy
type as $X \simeq_\Q X' \times K(\Q, n)$ for some simply connected
space $X'$.
\end{theorem}

\begin{proof}
The hypothesis  $SH_*(f) \neq 0$ translates into the following:
There is some $\alpha \in \pi_{n}(A)\otimes\Q$ whose image under
the composition $H_*(f)\circ h_\Q \colon \pi_{n}(A)\otimes\Q \to
H_{\text{n}}(X;\Q)$ is non-zero.  Therefore, $f_\#(\alpha) \in
\pi_{n}(X)\otimes\Q$ gives a non-zero spherical element
$h_\Q\big(f_\#(\alpha)\big) \in SH_{\text{n}}(X;\Q)$.  But since
$f$ is rationally cyclic, we have $f_\#(\alpha) \in
G_{n}(X)\otimes\Q$, with $h_\Q\big(f_\#(\alpha)\big) \not= 0$. By
a theorem of Oprea (see \cite[Lem.1.1]{Hal2}), this implies the
splitting $X \simeq_\Q X' \times S^{2n+1}$.
\end{proof}

Note that \thmref{thm:spherical on cohomology} can also be
interpreted as a  restriction on the homomorphism induced on
rational cohomology by a cyclic map:    A cyclic map into a space
not decomposable as above must have zero image in rational
spherical cohomology. We are unsure whether or not the
homomorphism induced on rational cohomology by a cyclic map may
contain indecomposable terms in its image without incurring the
splitting of \thmref{thm:spherical on cohomology}.

\begin{corollary}\label{cor:chi non-zero}
Suppose $X$ has finite dimensional rational homology and that $f
\colon A \to X$ is a rationally cyclic map. If $\chi(X) \not=0$,
then $SH_*(f) =0.$
\end{corollary}

\begin{proof}
From \thmref{thm:induced on even cohomology}, the image of
$H^*(f)$ in even degrees is decomposable, and hence non-spherical.
Thus $SH_*(f)$ is zero  in even degrees.
Since $\chi(X) \not=0$, the splitting of \thmref{thm:spherical on
cohomology} cannot occur with $n$ odd---recall that $K(\Q,2n+1)
\simeq_\Q S^{2n+1}$.  Therefore, by that result, $SH_*(f)$ is  zero
  in odd-degrees, as well.
\end{proof}

In particular, we see that under the hypotheses of \corref{cor:chi
non-zero}, the homomorphism induced by $f$ on rational
homology is zero if all homology of $H_{*}(A;\Q)$ is
spherical.  Thus, the corollary, or better, \thmref{thm:spherical
on cohomology}, gives a satisfactory generalization to cyclic maps
of Gottlieb's result concerning odd-degree Gottlieb groups.

From the above results, we can also extend the main result of
\cite{Lim}, at least in the simply connected case.  In the
original, very strong hypotheses are placed on $X$, namely that
both $H_*(X)$ and $\pi_*(X)$ are finitely generated. If,
in addition, we have $\chi(X) \not=0$ and $A$ a co-$H$-space, then
it is shown that a cyclic map $f \colon A \to X$ obtains $\Sigma
f$ an element of finite order in $[\Sigma A, \Sigma X]$.

Recall that a \emph{rational co-$H$-space} is a space $A$ such
that $A_\Q$ is a co-$H$-space. Recall that a \emph{(rationally)
elliptic space} is one whose rational homology and rational
homotopy are both finite dimensional.  A result of Halperin
\cite{Hal} states that the Euler characteristic of a rationally
elliptic space is non-negative. Further, that in the case of
positive Euler characteristic, that is, in the non-zero case, the
rational cohomology is zero in odd degrees. For brevity, we refer
to an elliptic space with positive Euler characteristic as an
\emph{$F_0$-space}. There are many interesting examples of such
spaces, including even dimensional spheres, complex projective
spaces, and ``maximal rank pair" homogeneous spaces $G/H$.
Products of $F_0$-spaces, and more generally many total spaces of
fibrations in which base and fibre are $F_0$-spaces, are again
$F_0$-spaces.

We can relax Lim's hypotheses as follows:

\begin{corollary}[cf.~{\cite[Th.5.2]{Lim}}]  Let $A$ be a
rational co-$H$-space with finite dimensional rational homology
and let $X$ be an $F_0$-space. Then any rationally cyclic map $f
\colon A \to X$ has rationally trivial suspension, $\Sigma f
\simeq_\Q
* \colon \Sigma A \to \Sigma X$.  If $A$ is assumed finite
dimensional, then the subgroup of cyclic maps $G(\Sigma A, \Sigma
X) \subseteq [\Sigma A, \Sigma X]$ is a finite group.
\end{corollary}

\begin{proof}
Since cup products in $\widetilde H^{*}(A;\Q)$ are trivial and
$H^*(X, \Q)$ is evenly graded we can apply \thmref{thm:induced on
even cohomology} to obtain $H^*(f) = 0 \colon \widetilde
H^{*}(X;\Q) \to \widetilde H^{*}(A;\Q)$. It follows that
$H_*(\Sigma f) = 0 \colon \widetilde H_*(\Sigma A;\Q) \to
\widetilde H_*(\Sigma X;\Q)$, and hence that $\Sigma f \simeq_\Q
*$. The assertion about finiteness follows from
\lemref{lem:finiteness}.
\end{proof}

\section{Examples}\label{sec:examples}

We give several examples of rationally cyclic maps. Our first
examples
illustrate that most  results from the previous section are sharp.
We then give some complete calculations of the set $G(A_\Q,X_\Q).$

Since our results tend towards restricting the possibilities for
rationally cyclic maps, we begin with an example that indicates
rich possibilities for such remain.  Of course, since any map into
an odd-dimensional sphere and any rational Gottlieb element are
rationally cyclic maps, we already have many examples.

\begin{example}
Consider the map $f \colon S^3 \times S^4 \to S^4$, obtained by
composing the quotient map $q \colon S^3 \times S^4 \to S^7$ with
the Hopf map $\eta \colon S^7 \to S^4$. It is well-known that $f$
is not rationally trivial, although it induces the trivial
homomorphism on both rational homotopy and rational cohomology.
Since $S^7$ is an $H$-space, $f$ itself, and not just its
rationalization, factors through an $H$-space. Actually, the
rationalization of $\eta$ is a rational Gottlieb element in
$G_7(S^4)\otimes\Q \cong \pi_7(S^4)\otimes\Q$, and so $f$ is a
rationally cyclic map by \thmref{thm:cyclic composition}. This
shows that we need $H^{\mathrm{odd}}(A;\Q) = 0$ in \thmref{thm:odd
rationally trivial}, even though the map may be trivial on
cohomology.  This example suggests many others of the form $A \to
S^n \to X$, with the second map a rational Gottlieb element.
\end{example}

\begin{example}
The image of $H^*(f)$ need not be decomposable in odd degrees, as
we now illustrate.  Take $A = S^2 \vee S^2 \cup_\alpha e^5$ with
the $5$-cell attached by the iterated Whitehead product $\alpha =
[\iota_1,[\iota_1,\iota_2]] \in \pi_4(S^2 \vee S^2)$. Then all cup
products in $H^*(A;\Q)$ are trivial, with $H^{5}(A;\Q)\cong\Q$
non-spherical.  Now the quotient map $q \colon A \to S^5$ is a map
into a rational $H$-space, and is therefore a rationally cyclic
map.  In degree $5$, $H^*(f)$ has indecomposable, non-spherical
image.
\end{example}

\begin{example}
Inclusion of the bottom cell $f \colon S^2 \to \C P^\infty$ is a
map into an $H$-space, and hence is a cyclic map. Evidently, we
have $H^{\mathrm{odd}}(S^2;\Q) = 0$ and yet $f$ is not rationally
trivial---it is not trivial on rational cohomology, in fact. This
shows that we need the hypothesis of finite rational category on
$X$ in \thmref{thm:odd rationally trivial}.  It also shows that
without $\cat0(X)<\infty$, then a cyclic map $f \colon A \to X$
may have $f_\#\otimes\Q$ non-zero in even degrees.  Indeed, here
we have $G_2(\C P^\infty) = \pi_2(\C P^\infty)$. This example
illustrates the splitting of \thmref{thm:spherical on cohomology}
for the case $n = 2$.
\end{example}

For some spaces $X$, including the $F_0$-spaces mentioned above,
the rational Gottlieb groups of $X$ coincide with the odd rational
homotopy groups of $X$.  In certain cases, this may generalize to
a bijection $G(A_\Q, X_\Q) \cong \bigoplus_{n \, \text{odd}}
H^n(A; \pi_n(X) \otimes \Q)$ for any space $A$.  We now give some
examples.

If $X$ is an odd dimensional sphere $S^{2n+1}$ then $X_\Q \simeq
K(\Q,2n+1)$ is an $H$-space and therefore $G(A_\Q,X_\Q) =
[A_\Q,X_\Q] = H^{2n+1}(A; \pi_{2n+1}(S^{2n+1}) \otimes \Q) =
H^{2n+1}(A;\Q)$. When $X = S^{2n}$ the same result holds but the
calculation is slightly more involved. The following example
includes the cases in which $X$ is an even dimensional sphere or a
complex projective space.

\begin{example}
Suppose $X$ has rational cohomology algebra a truncated polynomial
algebra on a single generator of even degree $2n$.  If $H^*(X;\Q)
\cong \Q[x]/(x^{k+1})$, then $G(A_\Q,X_\Q) \cong H^{2n(k+1)-1}(A;\Q)$.
We argue as follows: First, $X$ has minimal model $\Lambda(x,
y)$ with $|x| = 2n$, $|y| = 2n(k+1)-1$, $d_X(x) = 0$,
and $d_X(y) = x^{k+1}$. Suppose $f \colon A \to X$ is a rationally
cyclic map with minimal model $\mathcal{M}_f \colon \mathcal{M}_X
\to \mathcal{M}_A$ given on generators by $\mathcal{M}_f(x) = a$
and $\mathcal{M}_f(y) = b$, for $a, b \in \mathcal{M}_A$. Now
consider an affiliated map $F \colon A \times X \to X$ for $f$.
Then by taking into account the relative degrees of $x$ and $y$,
we can write its minimal model $\mathcal{M}_F \colon \mathcal{M}_X
\to \mathcal{M}_A\otimes\mathcal{M}_A$ on generators as
$\mathcal{M}_F(x) = x + a$ and $\mathcal{M}_F(y) = y + b + b_1 x +
\cdots + b_k x^k$, for $b_i \in \mathcal{M}_A$.  Since
$\mathcal{M}_F$ is a DG map, we can equate the terms that occur in
$d_A \mathcal{M}_F(y)$ with those that occur in $\mathcal{M}_F
d_X(y) = (x + a)^{k+1}$.  Thus we find that $d_A(b_k) = (k+1)a$.
Now define a DG homotopy $\Psi \colon \mathcal{M}_{X} \to
\mathcal{M}_{A}\otimes\Lambda(t,dt)$ on generators by setting
$$\Psi(x) = a(1-t) + \frac{1}{k+1} b_k dt$$
and
$$\Psi(y) = b - \frac{1}{k+1}b_ka^k + \frac{1}{k+1}b_ka^k
(1-t)^{k+1}.$$
This is easily checked to commute with the differentials, and so
defines a DG homotopy that starts at $\mathcal{M}_{f}$ and ends at
a map that is zero on the even-degree generator $x$.  So far, we
have argued that a cyclic map $f \colon A \to X$ has minimal model
given on generators by $\mathcal{M}_f(x) = 0$ and
$\mathcal{M}_f(y) = b$, up to homotopy.   For a map of this form,
$b$ must be a cycle in $\mathcal{M}_A$.  Furthermore, two such
maps are DG homotopic exactly when the cycles that $y$ is sent to
by each represent the same cohomology class of $H(\mathcal{M}_A)$.
This latter assertion is easily justified using the approach of
\cite{Ark-Lup}, for instance.
\end{example}

\begin{example}  Suppose $X= G/T$ where $G$ is a compact, connected Lie group and
$T$ is a maximal torus. Then  $G(A_\Q, X_\Q) \cong \bigoplus_{n \,
\text{odd}} H^n(A; \pi_n(X) \otimes \Q)$ for any space $A$.    For
in this case the minimal model of $X$ is of the form $\Lambda(x_1,
\ldots, x_n, y_1, \ldots, y_n)$ with $|x_i| = 2, d_X(x_i)= 0$ and
the $y_i$ of odd degree with $d_X(y_i)$ a decomposable polynomial
in the $x_i$. The rational cohomology of $X$ is evenly graded and
the cycles $x_i$ correspond to a space of generators.  See
\cite[Prop.15.16 and Ex.2,p.448]{F-H-T} for justification of
these assertions. By \thmref{thm:induced on even cohomology}, if
$f : A \to X$ is rationally cyclic, then $H^*(f) = 0$ since
$H^2(A; \Q)$ cannot contain non-zero decomposable terms. Thus the
minimal model $\mathcal{M}_f \colon \mathcal{M}_X \to
\mathcal{M}_A$ of $f$ satisfies $\mathcal{M}_f(x_i) = 0$ and
$\mathcal{M}_f(y_i) = \xi_i$ for cycles $\xi_i \in \mathcal{M}_A$
and the result follows as above.
\end{example}

\providecommand{\bysame}{\leavevmode\hbox
to3em{\hrulefill}\thinspace}
\providecommand{\MR}{\relax\ifhmode\unskip\space\fi MR }
\providecommand{\MRhref}[2]{%
  \href{http://www.ams.org/mathscinet-getitem?mr=#1}{#2}
} \providecommand{\href}[2]{#2}


\begin{thebibliography}{HMR75}

\bibitem[AL95]{Ark-Lup}
M.~Arkowitz and G.~Lupton, \emph{On finiteness of subgroups of
self-homotopy
  equivalences}, Cech Centennial Conference, Contemp. Math., vol. 181, Amer.
  Math. Soc., 1995, pp.~1--25.

\bibitem[FH82]{F-H}
Y.~F{\'e}lix and S.~Halperin, \emph{Rational {L}.-{S}. category
and its
  applications}, Trans.~Amer. Math.~Soc. \textbf{273} (1982), no.~1, 1--38.

\bibitem[FHT01]{F-H-T}
Y.~F{\'e}lix, S.~Halperin, and J.-C. Thomas, \emph{Rational
homotopy theory},
  Graduate Texts in Mathematics, vol. 205, Springer-Verlag, New York, 2001.

\bibitem[Got69]{Go1}
D.~H. Gottlieb, \emph{Evaluation subgroups of homotopy groups},
Amer. J. Math.
  \textbf{91} (1969), 729--756.

\bibitem[Hal77]{Hal}
S.~Halperin, \emph{Finiteness in the minimal models of
{S}ullivan},
  Trans.~Amer.~Math.~Soc. \textbf{230} (1977), 173--199.

\bibitem[Hal88]{Hal2}
\bysame, \emph{Torsion gaps in the homotopy of finite complexes},
Topology
  \textbf{27} (1988), 367--375.

\bibitem[HMR75]{H-M-R}
P.~Hilton, G.~Mislin, and J.~Roitberg, \emph{Localization of
nilpotent groups
  and spaces}, North-Holland Publishing Co., Amsterdam, 1975, North-Holland
  Mathematics Studies, No. 15, Notas de Matem\'atica, No. 55.

\bibitem[Lan75]{Lan}
G.~Lang, \emph{Localizations and evaluation subgroups}, Proc.
Amer. Math. Soc.
  \textbf{50} (1975), 489--494.

\bibitem[Lim82]{Lim}
K.~L. Lim, \emph{On cyclic maps}, J.~Austral.~Math.~Soc.
\textbf{32} (1982),
  349--357.

\bibitem[Opr95]{Oprea}
J.~Oprea, \emph{Gottlieb groups, group actions, fixed points and
rational
  homotopy}, Lecture Notes Series, vol.~29, Seoul National University Research
  Institute of Mathematics Global Analysis Research Center, Seoul, 1995.

\bibitem[Tan83]{Tan83}
D.~Tanr{\'e}, \emph{Homotopie rationnelle: mod\`eles de {C}hen,
{Q}uillen,
  {S}ullivan}, Lecture Notes in Mathematics, vol. 1025, Springer-Verlag,
  Berlin, 1983.

\bibitem[Var69]{Var}
K.~Varadarajan, \emph{Generalised {G}ottlieb groups}, J.~Indian
Math. Soc.
  \textbf{33} (1969), 141--164.

\end{thebibliography}
\end{document}